%% file: agt-4-44.tex
\documentclass{gtart_h}
\input agtout

\lognumber{44}
\volumenumber{4}
\volumeyear{2004}
\papernumber{44}
\pagenumbers{1041}{1044}
\received{14 June 2004} 
\revised{20 August 2004}
\accepted{7 September 2004}
\published{3 November 2004}

\usepackage{amsmath,amssymb,amscd}

\def\Addressesrr{\bigskip
{\small \parskip 0pt \leftskip 0pt \rightskip 0pt plus 1fil \def\\{\par}
\sl\theaddress\par
\medskip
\rm Email:\stdspace\tt\theemail\hfill\rm Received:\qua\receiveddate 
\quad Revised:\qua\reviseddate\par}}

\newtheorem{thm}{Theorem}[section]

\newtheorem{prop}{Proposition}[section]

\theoremstyle{definition}
\newtheorem{defn}{Definition}[section]

\begin{document}
\title{Mp-small summands increase knot width}
\author{Jacob Hendricks}
\address{Department of Mathematics, University of 
Arkansas, Fayetteville, AR 72701, USA}
\email{jghendr@uark.edu}

\begin{abstract}
{Scharlemann and Schultens have shown that for any pair of knots $K_1$
and $K_2$, $w(K_1 \# K_2) \geq max\{w(K_1), w(K_2)\}$.
Scharlemann and Thompson have given a scheme for possible examples where
equality holds.  Using results of Scharlemann--Schultens,
Rieck--Sedgwick and Thompson, it is shown that for $K= \#_{i=1}^n
K_i$ a connected sum of mp-small knots and $K^{\prime}$ any
non-trivial knot, $w(K\#K^{\prime})>w(K)$. }
\end{abstract}

\asciiabstract {Scharlemann and Schultens have shown that for any pair of
knots K_1 and K_2, w(K_1 \# K_2) >= max{w(K_1),w(K_2)}.
Scharlemann and Thompson have given a scheme for possible examples where
equality holds.  Using results of Scharlemann-Schultens,
Rieck-Sedgwick and Thompson, it is shown that for K= \#_{i=1}^n
K_i a connected sum of mp-small knots and K' any
non-trivial knot, w(K\#K')>w(K).}


\makeshorttitle

\section{Introductory remarks}
Thin position, introduced by David Gabai \cite{g}, has applications in
solving difficult problems, yet there remains much to understand
about the width of knots.   Current understanding of the behavior of
width under connect sum is incomplete, but has some form.
For example, given two knots $K_1$ and $K_2$, it is easily seen that
$w(K_1 \# K_2) \leq w(K_1) + w(K_2) -2$;
simply stack the knots in thin position, perform the connect
sum, and calculate the upperbound.
 Also, Rieck and Sedgwick \cite{rs} have shown that if $K_1$ and $K_2$ are mp-small
 (see definition \ref{11}) then $w(K_1 \# K_2) = w(K_1) + w(K_2) -
2$.
By means of a fascinating method (see \cite{ss}), it has also been shown that for
any two knots $K_1$ and $K_2$, $w(K_1 \# K_2) \geq \mbox{max}\{w(K_1),
w(K_2)\}$.  It is natural to wonder if there exists a pair of knots
such that $w(K_1 \# K_2) = w(K_1)$, for this would be a peculiar
property of width.  In fact, a class of possible examples of equality is presented in
\cite{st}, where it is argued, though not proven, that for any
knot $K_2$ there is a (typically quite complicated) knot $K_1$,
apparently determined by the bridge number of $K_2$, so that $w(K_1 \# K_2) =
w(K_1)$.  In contrast, we show here that if $K_1$ is the connect sum of mp-small knots,
then for any knot $K_2$, $w(K_1 \# K_2) > w(K_1)$.

I would like to thank the referees for many helpful remarks and Yo'av
Rieck for many helpful conversations.

\section{Preliminaries}

For the sake of brevity, familiarity with \textit{width, thin/thick
levels, thin position, bridge position, swallow-follow torus, satellite knot}
as well as the following definitions and theorems will be
assumed. One can find definitions of the
italicized words in many sources, for example \cite{hk} and \cite{ss}.

\begin{defn}
A knot $K \subset S^3$ is called meridionally planar small (mp-small) if the
only incompressible meridional surface in its complement is a boundary-parallel
annulus.
\label{11}
\end{defn}

In \cite{t}, Thompson proved:

\begin{thm}[Thompson]
If a knot $K$ in thin position countains a thin level, then it is
not mp-small.
\label{1}
\end{thm}

Thus, for an mp-small knot $K$, thin position must equal bridge
position. In \cite{rs} (Theorem 4.1), a converse of Theorem \ref{1}
is given in the case that $K$ is the connect sum of two non-trivial knots.

\begin{thm}[Rieck--Sedgwick]
Let $K$ be a connected sum of non-trivial knots, $K = K_1 \# K_2$.
Then any thin position for $K$ is not bridge position for K.
\label{3}
\end{thm}

The following was also shown in \cite{rs}:

\begin{thm}[Rieck--Sedgwick]
Let $K= \#_{i=1}^n K_i$ be a connected sum of mp-small knots. If K is
in thin position, then there is an ordering of the summands $K_{i_1},
K_{i_2},\ldots ,K_{i_n}$ and a collection of leveled decomposing annuli
$A_{i_1}, A_{i_2},\ldots ,A_{i_n-1}$ so that the thin levels of the
presentation are precisely the annuli $\{ A_{i_j} \}$ occurring in
order, where the annulus $A_{i_j}$ separates the connected sum
$K_{i_1} \# K_{i_2} \#\ldots \# K_{i_j}$ from the connected sum $K_{i_j+1}
\#\ldots \# K_{i_n}$.
\label{rs2}
\end{thm}

Scharlemann--Schultens \cite{ss} describe a method for reimbedding a knot $K$ so that a height
function on $K$ would be preserved.

\begin{thm}[Scharlemann--Schultens]
Suppose $p: S^3 \rightarrow \mathbb{R}$ is the standard height
function. Let $K \subset S^3$ be a knot that lies inside a
standard unknotted torus $H \subset S^3$. Let $f: H \rightarrow S^3$ be a
possibly knotted embedding of $H$ in $S^3$.  Then there is a
reimbedding $f^{\prime}: H \rightarrow S^3$ so that

\begin{enumerate}
\item  $pf^{\prime} = pf$, i.e. the reimbedding preserves height
\item  $S^3 - f^{\prime}(H)$ is a solid torus (so $f^{\prime}(H)$ is an
unknotted solid torus) and
\item  $f(K^{\prime})$ is isotopic to $K$ in $S^3$.
\end{enumerate}
\label{2}
\end{thm}

This is a specific case of what is shown in Corollary 5.4 of
\cite{ss}. In this paper will apply the theorem in this setting:
We will be given a connected sum $K\#L$ embedded in $S^3$ via a
knotted embedding $f: H \rightarrow S^3$.  $K$ has here been placed in $H$
as a wrapping number one knot and the embedding $f$ is known to
take the core of $H$ to the knot $L$.  (The boundary of $f(H)$ is
commonly called a "swallow-follow" torus for $K\#L$, for it
swallows $K$ and follows $L$.)

\section{When summands are mp-small}

In \cite{ss}, it is shown that $w(K_1\#K_2)\geq w(K_i)$ for $i=1,2$.
One may ask if there is a pair of knots $K_1$ and $K_2$ such that equality
holds (eg. \cite{st}).

It is shown in \cite{rs} that for $K_1$ and $K_2$ both mp-small
knots, $w(K_1\#K_2) = w(K_1)+w(K_2)-2$.  In particular, when both
are non-trivial knots, $w(K_1\#K_2)$ $>w(K_i)$, $i=1,2$.  A natural
question is whether this inequality remains true if only one of
the knots is mp-small.

\begin{prop}
If $K_1$ is an mp-small knot, then for any non-trivial knot $K_2$,
$w(K_1\#K_2)>w(K_1)$.
\label{1case}
\end{prop}

\begin{proof}

Since $K_1$ is mp-small, Theorem \ref{1} implies that thin position of $K_1$ equals
bridge position. Then, take $K_1\#K_2$ to be a satellite knot with
companion knot $K_2$ via the \textit{swallow-follow torus}.  Put $K_1\#K_2$ in
 thin position;
notice by Theorem \ref{3}, $K_1\#K_2$ cannot be in bridge
position, so there exists a thin level.
By Theorem \ref{2}, there exists a height preserving reimbedding that changes
the knot $K_1\#K_2$ to the knot $K_1$. Since $K_1\#K_2$ has a thin level, so must
a height preserving reimbedding of $K_1\#K_2$.
So, this reimbedding yields an embedding of $K_1$ that has a thin level.
 (i.e. an embedding of $K_1$ that is not in bridge
position.) It follows that this embedding of $K_1$ cannot be
in thin position, so its width as embedded, and hence the width of $K_1\#K_2$,
 is greater than the minimal width of $K_1$.
\end{proof}

One can use the method for showing proposition \ref{1case} to prove a more general
statement; that is, $K_1$ can be taken to be the connect sum of
mp-small knots.

\begin{thm}
Let $K= \#_{i=1}^n K_i$ be a connected sum of mp-small knots and let $K^{\prime}$
be any non-trivial knot.  Then, $w(K\#K^{\prime})>w(K)$.
\end{thm}

\begin{proof}

Suppose $w(K\#K^{\prime})=w(K)$; we will show that $K^{\prime}$ must be the unknot.
Now, take $K^{\prime}$ to be the companion knot of the satellite knot given by
$K\#K^{\prime}$ with $H$ a torus that swallows $K$ and follows $K^{\prime}$.
Put $K\#K^{\prime}$ in thin position. By Theorem \ref{2}, there exists a reimbedding
$f$ such that $f$ preserves a height function on $K\#K^{\prime}$ and
$f(K\#K^{\prime})$ is isotopic to $K$ in $S^3$.  We are assuming
$w(K\#K^{\prime})=w(K)$, thus $f(K\#K^{\prime})$ is a thin
presentation of $K$; therefore, by Theorem \ref{rs2}, every thin level of $K$
is one of the $n-1$ decomposing annuli, which decompose the
connect sum $\#_{i=1}^n K_i$ into mp-small knots.  Since each $K_i$
is mp-small, by Theorem \ref{1}, thin position of each $K_i$ is bridge
position of each $K_i$.
Now, since $f$ is a height preserving reimbedding, $K\#K^{\prime}$ must
have $n-1$ decomposing annuli that constitute all of the thin
levels; hence, $K\#K^{\prime} = \#_{i=1}^n K^{\prime}_i$.  Note that each
 of the $n$ components $K^{\prime}_i$
is in bridge position.  This must also be thin position; for if
not, one could thin $K\#K^{\prime}$ by thinning the summand for which
bridge does not equal thin, but we have put $K\#K^{\prime}$ in thin
position.  It follows from Theorem \ref{3} that each $K^{\prime}_i$ must be
prime.  Hence, we have $K\#K^{\prime} = \#_{i=1}^n K^{\prime}_i$ where each
$K^{\prime}_i$ is prime, and $K=\#_{i=1}^n K_i$ where each $K_i$ is prime.
Since knot factorizations are unique, $K^{\prime}$ must be the unknot.
\end{proof}

\Addressesrr

\end{document}

%% file: agtout.tex
\def\ifplaintex{\expandafter\ifx\csname documentclass\endcsname\relax}

\def\gtp{{\mathsurround=0pt\it $\cal G\mskip-2mu$eometry \&\ 
$\cal T\!\!$opology $\cal P\!$ublications}}  

\def\recd{{\small Received:\qua\receiveddate\ifx\reviseddate\relax
\else\qquad Revised:\qua\reviseddate\fi\par}} 

\def\lognumber#1{\def\thelognumber{#1}}
\def\volumenumber#1{\def\thevolumenumber{#1}}
\def\volumeyear#1{\def\thevolumeyear{#1}}
\def\papernumber#1{\def\thepapernumber{#1}}
\def\pagenumbers#1#2{\def\startpage{#1}\def\finishpage{#2}}
\def\published#1{\def\publishdate{#1}}

\def\received#1{\def\receiveddate{#1}}
\def\revised#1{\def\reviseddate{#1}}
\def\accepted#1{\def\accepteddate{#1}}

\long\def\asciiabstract#1{\long\def\theasciiabstract{#1}}

\let\\\par\let\thelognumber\relax\let\thevolumenumber\relax
\let\thepapernumber\relax\let\thevolumeyear\relax\let\startpage\relax
\let\finishpage\relax\let\publishdate\relax\let\receiveddate\relax
\let\reviseddate\relax\let\accepteddate\relax\let\theasciititle\relax
\let\theasciiauthors\relax
\let\theasciiabstract\relax

\let\theasciiemail\relax

\ifplaintex
\font\logobig=cmssbx10 scaled 3836
\font\logomed=cmssbx10 scaled 2557
\else
\font\logobig=cmssbx10 scaled 4200
\font\logomed=cmssbx10 scaled 2800
\fi

\long\def\makeagttitle{   
\count0=\startpage
\agt\hfill      
\hbox to 45truept{\vbox to 0pt{\vglue -13truept{\logomed A\kern -.37em{\logobig 
T}\kern -.38em G}\vss}\hss}
\break
{\small Volume \thevolumenumber\ (\thevolumeyear)
\startpage--\finishpage\nl
Published: \publishdate}

\vglue .25truein

{\parskip=0pt\leftskip 0pt plus
1fil\def\\{\par\smallskip}{\Large\bf\thetitle}\par\medskip} \vglue
0.05truein

{\parskip=0pt\leftskip 0pt plus 1fil\def\\{\par}{\sc\theauthors}
\par\medskip}%
 
\vglue 0.03truein 

{\small\leftskip 25truept\rightskip 25truept{\bf Abstract}\stdspace\theabstract

{\bf AMS Classification}\stdspace\theprimaryclass
\ifx\thesecondaryclass\relax\else; \thesecondaryclass\fi\par
{\bf Keywords}\stdspace \thekeywords\par}\vglue 7truept

}   

\ifplaintex
\font\phead=cmsl9 scaled 950
\font\pnum=cmbx10 scaled 913
\font\pfoot=cmsl9 scaled 950
\headline{\vbox to 0pt{\vskip -4.5mm\line{\small\phead\ifnum
\count0=\startpage ISSN 1472-2739 (on-line) 1472-2747 (printed)
\hfill {\pnum\folio}\else\ifodd\count0\def\\{ }%
\ifx\theshorttitle\relax\thetitle\else\theshorttitle\fi\hfill{\pnum\folio}
\else\def\\{ and }{\pnum\folio}\hfill\ifx\theshortauthors\relax\theauthors
\else\theshortauthors\fi\fi\fi}\vss}}
\footline{\vbox to 0pt{\vglue 0mm\line{\small\pfoot\ifnum\count0=\startpage
\copyright\ \gtp\hfill\else
\agt, Volume \thevolumenumber\ (\thevolumeyear)\hfill\fi}\vss}}
\else
\font=cmsl9 scaled 1050
\font\lnum=cmbx10 
\font=cmsl9 scaled 1050
\makeatletter
\def\@oddhead{{\small\lhead\ifnum\count0=\startpage ISSN 1472-2739 
(on-line) 1472-2747 (printed)\hfill {\lnum\number\count0}\else\ifodd\count0
\def\\{ }\ifx\theshorttitle\relax \thetitle \else\theshorttitle\fi\hfill
{\lnum\number\count0}\else\def\\{ and }{\lnum\number\count0}
\hfill\ifx\theshortauthors\relax 
\theauthors\else\theshortauthors\fi\fi\fi}}\def\@evenhead{\@oddhead}
\def\@oddfoot{\small\lfoot\ifnum\count0=\startpage\copyright\ \gtp\hfill\else
\agt, Volume \thevolumenumber\ (\thevolumeyear)\hfill\fi}
\def\@evenfoot{\@oddfoot}
\makeatother
\fi
\let\maketitlepage\makeagttitle
\let\makeshorttitle\maketitlepage
\let\maketitle\maketitlepage

\newwrite\gtoutfile
\long\gdef\makeheadfile{  
{\def\\{, }\def\s{ }
\immediate\openout\gtoutfile head.xxx
\immediate\write\gtoutfile{Proxy-for: \ifx\theasciiauthors\relax
\theauthors\else\theasciiauthors\fi\s<\ifx\theasciiemail\relax\theemail\else\theasciiemail\fi>}
\immediate\write\gtoutfile{\noexpand\\}
\immediate\write\gtoutfile{Authors: \ifx\theasciiauthors\relax
\theauthors\else\theasciiauthors\fi}
{\def\\{ }\immediate\write\gtoutfile{Title: \ifx\theasciititle\relax
\thetitle\else\theasciititle\fi}}
\immediate\write\gtoutfile{Subj-class: GT or SG, GR etc}
\immediate\write\gtoutfile{MSC-class: \theprimaryclass\ifx\thesecondaryclass\relax\else, \thesecondaryclass\fi}
\immediate\write\gtoutfile{Journal-ref: Algebr. Geom. Topol. \thevolumenumber\s
(\thevolumeyear) \startpage-\finishpage}
\immediate\write\gtoutfile{Comments: Published by Algebraic and
Geometric Topology at}
\immediate\write\gtoutfile{\s\s\s  http://www.maths.warwick.ac.uk/agt/AGTVol\thevolumenumber/agt-\thevolumenumber-\thepapernumber.abs.html}
\immediate\write\gtoutfile{\noexpand\\}
\immediate\write\gtoutfile{}
\ifx\theasciiabstract\relax
\immediate\write\gtoutfile{\theabstract}\else
\immediate\write\gtoutfile{\theasciiabstract}\fi
\immediate\write\gtoutfile{}
\immediate\write\gtoutfile{\noexpand\\}
\immediate\write\gtoutfile{}
\immediate\closeout\gtoutfile}}  

\def\maketitlepage{\makeagttitle\makeheadfile}
\let\makeshorttitle\maketitlepage
\let\maketitle\maketitlepage

%% file: agt-4-44.bbl
\vspace{-0.5cm}

\begin{thebibliography}

\itemsep 1pt plus 1pt
\bibitem{g} \textbf{D Gabai}, {\em Foliations and the topology of
3-manifolds, III}, J. Differential Geom. 26 (1987) 479-536
\MR{0910017}

\bibitem{hk} \textbf{D\,J Heath}, \textbf{T Kobayashi},
{\em{Essential Tangle Decomposition from Thin Position of a Link}}, 
Pacific J. Math. 179 (1997)  101--117
\MR{1452527}

\bibitem{rs} \textbf{Y Rieck}, \textbf{E Sedgwick}, {\em{Thin position
 for a connected sum of small knots}}, \agtref2{2002}{14}{297}{309}
\MR{1917054}

\bibitem{ss} \textbf{M Scharlemann}, {\bf J Schultens},
 {\em 3-manifolds with planar presentations and the width
of satellite knots}, \arxiv{math.GT/0304271}

\bibitem{st} \textbf{M Scharlemann}, {\bf A Thompson}, {\em On the
 additivity of knot width}, from: ``Proceedings of the Casson Fest'',
 (Cameron Gordon and Yoav Rieck, editors), \gtmref7{2004}5{135}{144}

\bibitem{t} \textbf{A Thompson}, {\em Thin position and bridge number
for knots in the 3-sphere}, Topology 36 (1997) 505--507
\MR{1415602}

\end{thebibliography}
